\documentclass[12pt]{article}
\UseRawInputEncoding
\usepackage{amsmath}
\usepackage{amssymb}
\usepackage{amsthm}
\usepackage{amsfonts}
\usepackage{amscd}
\usepackage[mathscr]{eucal}
\newcommand{\Z} {{\mathbb  Z}}
\newcommand{\Q}{{\mathbb  Q}}

\begin{document}
\parindent  25pt
\textwidth  15cm    \textheight  23cm \evensidemargin -0.06cm
\oddsidemargin -0.01cm

\title{ { On Some Diophantine Parameters of the Cyclic Torsion Subgroups
of Odd Order of Elliptic Curves over  $ \Q $  }}
\author{\mbox{}
{Derong Qiu}
\thanks{ \quad E-mail:
derong@mail.cnu.edu.cn, \ derongqiu@gmail.com } \\
(School of Mathematical Sciences,
 Capital Normal University, \\
Beijing 100048, P.R.China )  }

\date{}
\maketitle
\parindent  24pt
\parskip  0pt

\par   \vskip 0.3cm

{\bf Abstract} \quad For elliptic curves $ E $ over the rational number field $ \Q, $
based on Mazur's classical result, we give in this paper some explicit Diophantine parameters
for their cyclic torsion subgroups $ E(\Q)_{\text{tors}} $ of odd order. This work is a
continuation of the ones in ([Q1], [QZ1$\sim $2]). Such explicit form of parameters might be
useful in calculating torsion subgroups of elliptic curves over higher
degree number fields.

\par  \vskip  0.2cm


\par     \vskip    1cm

 \hspace{-0.6cm}{\bf 1 \ \ Introduction and Main Results}

      \par \vskip 0.6 cm

Let $ E $ be an elliptic curves over the rational number field $
\Q. $  From Mordell-Weil theorem, the set $ E (\Q)$ of
rational points is a finitely generated abelian group. By a deep
theorem of Mazur, the torsion subgroup $ E(\Q)_{\text{tors}} $ of $
E (\Q )$ is isomorphic to one of the following fifteen groups
(see [M1],[M2],[S]):
$$ (\text{Cyclic} N-\text{types}) \quad \quad \quad \quad \Z/N\Z \qquad
(1\le N\le 10 \quad \text{or} \quad N = 12); $$
$$ (\text{Non-cyclic} \ (2, \ 2N)-\text{types}) \quad \quad \quad  \quad \Z/2\Z
\times \Z /2N\Z \qquad (1\le N\le 4 ). $$
\par   \vskip 0.1cm

In 1996, a kind of explicit parameters of the non-cyclic $ (2, 2N )-$torsion types of
$ E(\Q)_{\text{tors}} $ was given by K.Ono ([O]), and in 1999 another parameters corresponding to
cyclic $ N-$torsion types for even $ N $  were given by Qiu and Zhang ([Q1],[QZ1$\sim $2]).
Such explicit parameters may be conveniently used in further study of the structure of rational
points of elliptic curves. And recently, similar results have been applied to determine the torsion
structure of elliptic curves over quadratic number fields, multi-quadratic number fields and the
$ \Z_{p}-$extensions of $ \Q $ (see [F1$\sim $3], [Kw] and [Q1],[QZ3]).
\par \vskip 0.2cm

In this paper, we continue to study the problem of explicit parameters associated to
rational torsion groups of elliptic curves, that is,
we consider the case when $ E(\Q)_{\text{tors}} $ is of the cyclic $
N -$types for odd $ N, $ \ and give a kind of explicit parameters for them.
\par \vskip 0.2cm

For the elliptic curve $ E $ over $ \Q, $ up to $ \Q-$isomorphism, $ E $ has a model
$$ E = E_{(a,b)}: \  y^{2}= x^{3} + ax + b ,  \quad  \quad \quad  a, \
b \in \Z . $$ Let $ p= (x, y ) \in E( \Q )_{tors} $ be a $ \Q- $
torsion point of $ E, $ then by Lutz-Nagell theorem (see [Kn] or
[S]), we always have $ x, \ y \in \Z . $ Let $ f(x) = x^{3} + ax
+ b . $  It is easy to see that $ E(\Q)_{\text{tors}} $ has no
non-trivial 2-torsion points (i.e. $ E(\Q)[2] = \{O\} $) if
and only if $ f(x) $ has no roots in $ \Q. $ \ For an abelian
group $ A  $  and a positive rational integer $ n , $ we denote $
A[n] = \{ a \in A: \ n a = 0 \}. $ \ Throughout this paper, $ O
\in E ( \Q ) $ is the zero element.
\par \vskip 0.2cm

Now we state our main results.
\par \vskip 0.2cm

{\bf Theorem 1.}\quad  Let $ E = E_{(a,b)}: \  y^{2}= x^{3} + ax
+ b \ $ be an elliptic curve with $  a, b \in \Z , $ \ and $ E
(\Q)[2] = \{O \} . $ Then the $ \Q-$rational torsion subgroup
$ E(\Q)_{\text{tors}} $ is parameterized as follows:  \\
(I) $ \ E(\Q)_{\text{tors}} \supseteq {\Z} /3{\Z} $ \quad if
and only if
$$ a=6mn - 27 n^{4}, \qquad   b = m^{2} - 18mn^{3} + 54 n^{6}, $$
 where  $ m, \ n \in {\Z }  $ and $ m \neq 0. $  \\
(II) $ \ E(\Q)_{\text{tors}} = {\Z}/ 9{\Z} $ \quad if and only if
$$ a = 6mn - 27 n^{4}, \qquad  b = m^{2} - 18mn^{3} + 54 n^{6},
\quad  m, \ n \in {\Z}, \ m \neq 0, $$ and $ m, n $ satisfying
the following conditions:  \\
There exist $ u, v, w \in \Z $ and $ v \neq 0 $ such that
$$ 2mw + 6mn = u^{2} + u v, \quad  w^{2} = 9 n^{2} + 2u + v, $$
$$  u \left[ ( 2 u + v )^{2} + 2uv  \right] + 4 ( m + 3 n u ) \left[
m + 3 n ( u + v ) \right] = 0. $$
(III) $ \ E(\Q)_{\text{tors}} = {\Z}  / 5{\Z} $ \quad if and only if
$$ a = -27 \left[( m^{2}-n^{2} + 6mn)^{2} - 20 m^{2} n^{2} \right], $$
$$ b = 54 ( m^{2}+ n^{2} ) \left[ ( m^{2}-n^{2} + 9mn )^{2} - 5 m^{2} n^{2} \right], $$
where $ m, n \in \Z $ and  $ m n \neq 0. $  \\
(IV) $ \ E(\Q)_{\text{tors}} = {\Z} / 7{\Z} $ \quad if and only if
 $$ a = -3 \lambda ^{2} + 4w ( u^{2} - v^{2} - w ), $$
$$ b = 8w^{2} ( u^{2} + v^{2} - \lambda ) + 2\lambda ( \lambda ^{2} +
2w ( v^{2} - u^{2} )), $$
$$ 3 \lambda + 2w = ( u - v )^{2}, \qquad  w^{2} + uvw - uv^{3} = 0, $$
where $ \lambda , u , v, w \in \Z $ and $ \ uvw \neq 0. $  \\
(V) In all other cases, $ E(\Q)_{\text{tors}} = \{ O \}. $
                 \par \vskip 0.2cm

Furthermore, a generator $ P_{ n } $  of the cyclic group $ E(\Q)_{\text{tors}} $ of
order $ n $ is obtained. In each case as above,  $ P_{ n } $ and $ 2P_{ n }$ are as follows:
\par \vskip 0.24cm

(I) $ P_{3} = ( 3 n^{2}, \ m ) ; \quad \quad  2P_{3} = ( 3 n^{2}, \ - m). $ \\
(II) $ P_{9} = (3 n^{2} + u, \  u w - m) ;  \quad \quad 2P_{9}
= (3 n^{2} + u + v, \ ( u + v ) w - m ). $  \\
(III)  $ P_{5} = ( 3 ( m^{2} + n^{2} ) + 18 mn, \ 108 m n^{2});
 \quad  \quad 2P_{5} = (3 (m^{2} + n^{2} ) - 18 mn, \ 108 m^{2} n ). $  \\
(IV) $ P_{7} = ( \lambda + 2w, \ 4 uw); \quad \quad 2P_{7} = (\lambda - 2w, \ -4 vw). $  \\
(V) $ P_{1} = O. $
\par  \vskip 0.2 cm

{\bf Remark 1.1.} \ It is obvious that the variable \ $ v $ \
in the equations of case ( II ) is superfluous, so does one of the
two variables \ $ \lambda $ \ and \ $ w $ \ in the equations of
case (IV). Since omitting them will make the expressions of the
corresponding equations become more complicated, we would rather
not change them.

{\bf Remark 1.2.} \ This work is a continuation and completion
of our ones in ([Q1], [QZ1$\sim $2]) many years ago. For a given family
of elliptic curves, there may be different ways to parameter them
according to their models and invariants. Such explicit form of parameters might be
useful in calculating torsion subgroups of elliptic curves over higher
degree number fields (see e.g. [F1$\sim $3],[Kw],[QZ3]).
\par  \vskip 0.3 cm

\hspace{-0.6cm}{\bf  2 \ \ Proof of the Theorem}
\par \vskip 0.2 cm

{ \bf Proof of Theorem 1.} \ By Lutz-Nagell Theorem (see [S]) we know that any
point $ P =( x, y ) \in E(\Q)_{\text{tors}} $ is an integer point, i.e., $ x, y \in \Z. $
Also, by our assumption, $ y \neq 0. $
\par \vskip 0.2 cm

(I) \ If $ \ E(\Q)_{\text{tors}} \supseteq \Z / 3\Z. $ \ Then there
exists a point $ P_{3} = ( x, \ y) \in E(\Q) $ such that $ 3P_{3} = O $ and
$ P_{3} \neq O . $ So $ 2 P_{3} = - P_{3} $ and $ x( 2 P_{3} ) = x( -P_{3} ) = x ( P_{3} ) = x. $
By the duplication formula (see [S], p.59)
$$ x ( 2P ) = \frac{ x^{4} - 2a x^{2} - 8 b x + a^{2}}{ 4 x^{3} + 4 a x + 4 b }.  \eqno(2.1) $$
 Since $ x ( 2 P_{3}) = x, $ by (2.1) we get
 $$ 3 x^{4} + 6 a x^{2} + 12 b x = a^{2} .   \eqno(2.2) $$
Thus $ 3 | a , $ so $ a = 3 a_{0} $ for some $ a_{0} \in \Z . $
From (2.2), $ ( x^{2} + 3 a_{0} )^{2} = 4 ( 3 a_{0}^{2} - b x ) .
$ Hence  $ 3 a_{0}^{2} - b x = c^{2} $ and $ x^{2} + 3 a_{0} = 2 c
$ for some $ c \in \Z . $ So $ a = 3 a_{0} = 2 c - x^{2} $ and $ 3
| ( 2 c - x^{2}) . $ Thus $ ( 2 c - x^{2})^{2} = 9 a_{0}^{2} = 3 (
b x + c^{2} ) , $ i.e., $ c^{2} - 4 c x^{2} + x^{4} = 3 b x . $ On
the other hand, by definition, $ b = y^{2} - x^{3} - a x = y^{2} -
x^{3} - ( 2 c - x^{2} ) x = y^{2} - 2 c x . $ Therefore $ c^{2} -
4 c x^{2} + x^{4} = 3 x ( y^{2} - 2 c x ) , $ i.e., $$ ( c + x^{2}
)^{2} = 3 x y^{2} .  \eqno( 2.3 ) $$
So $ 3x = n_{0}^{2} $ for some $ n_{0} \in \Z, $ i.e., $ 3 \mid n_{0}^{2}, $
hence $ 3 \mid n_{0}, $ and $ n_{0} = 3 n $ for some $ n \in \Z . $ Therefore
$ x = 3 n^{2}. $ Substituting this $ x $ into (2.3), we get $ c = 3 n y - 9
n^{4} . $ Let $ y = m \in \Z \setminus \{ 0 \} . $ Then
$$ a = 2 c - x^{2} = 2 ( 3 n m - 9 n^{4} ) - ( 3 n^{2} )^{2} = 6 m n -
27 n^{4} , $$    $$ b = y^{2} - 2 c x = m ^{2} - 2 ( 3 n m - 9
n^{4} ) ( 3 n^{2} ) = m^{2} - 18 m n^{3} + 54 n^{6} . $$
Conversely, if the conditions on $ a, b $ in (I) hold, then
from the above procedure, it is easy to verify that the point $
 P_{3} = (3 n^{2},  m) \in E (\Q)_{\text{tors}} $ is of order 3.
 Also $ 2 P_{3} = - P_{3} = ( 3 n^{2}, -m ). $
 \par \vskip 0.2cm

(II) \ If $ E(\Q)_{\text{tors}} = \Z / 9 \Z. $ Then $ E(\Q)_{\text{tors}}
\supseteq \Z / 3 \Z, $ \ and there exists a point $ P_{9} = ( x, \ y ) \in
E (\Q)_{\text{tors}} $ of order 9. So by case (I),
$$ a = 6 m n - 27 n^{4}, \quad  \quad b = m^{2} - 18 m n^{3}
+ 54 n^{6}  $$ for some $ m, n \in \Z $ and $ m \neq 0; $ And $
\ E( \Q ) [3] = \{O, P_{3}, 2 P_{3}\}, \ $ where $ P_{3} = ( 3 n^{2}, \ m ). $ \
Let $ P_{3}^{ \prime } = ( x_{3}, \ y_{3}) = 3 P_{9}, $ \ then
$ P_{3}^{ \prime } \in E( \Q ) [3] $ is of order 3. So $ x_{3} = x( P_{3}) = x( 2 P_{3}) = 3 n^{2} $
and $ y_{3} = m \varepsilon $ with $ \varepsilon \in \{ 1, -1 \}. $ \
By the addition law (see [S], pp.58$\sim$59),
$$ x_{3} = x ( P_{9} + 2 P_{9} ) = \left( \frac{y_{2} - y }{ x_{2} - x } \right)^{2}
 - x - x_{2},   \qquad \qquad \text{and} $$
 $$ y_{3} = y( P_{9} + 2P_{9}) = - \left( \frac{y_{2} - y }{
x_{2} - x } \right) x_{3} - \frac{ x_{2} y - x y_{2}}{ x_{2} - x }
= \left( \frac{y_{2} - y }{ x_{2} - x } \right) ( x - x_{3} ) - y, $$
here $ 2P_{9} = ( x_{2} , \ y_{2}). $
 Obviously, $ x_{2} \neq x . $ \ Let $ w = \frac{y_{2} - y }{ x_{2} - x} . $
Since $ x, y, x_{2}, y_{2}, x_{3} \in \Z, $ we have $ w \in
\Z , $ and then
$$ \left \{
   \begin{array}{l}
  x + x_{2} + 3 n^{2} = w^{2} , \\
  y = w ( x - 3 n^{2} ) - m \varepsilon , \\
  y_{2}= w ( x_{2} - x ) + y = w ( x_{2} - 3 n^{2} ) - m \varepsilon .
  \end{array}
  \right.   \eqno(2.4) $$
Thus $ y + y_{2} = w ( w^{2} - 9 n^{2} ) - 2 m \varepsilon . $
Therefore by definition, $$ \begin{array}{l} ( x_{2} - x ) ( x^{2}
+ x x_{2} + x_{2}^{2} + a ) = ( x_{2}^{3} + a x_{2} + b  ) - (
x^{3} + a x + b) \\
= y_{2}^{2} - y^{2} = ( y_{2} - y ) ( y_{2} + y ) \\
 = w( x_{2} - x ) ( w ( w^{2} - 9 n^{2} ) - 2 m \varepsilon ).
 \end{array} $$  Since $ x_{2} - x \neq 0, $ we get
 $$ ( x +  x_{2})^{2} - x x_{2} + a = w( w ( w^{2} - 9 n^{2} ) - 2 m \varepsilon ),
 \qquad \text{ so } $$
$$ \begin{array}{l} x x_{2} = ( w^{2} - 3 n^{2} )^{2} + ( 6 m n
- 27 n^{4}) - w( w ( w^{2} - 9 n^{2} ) - 2 m \varepsilon ) \\
= 6 m n - 18 n^{4} + 3 n^{2} w^{2} + 2m w \varepsilon .
\end{array}  \eqno(2.5) $$
Since $ x $ and $ x_{2} $ are the two distinct integer roots of
the polynomial $$ h(T) = T^{2} - ( w^{2} - 3 n^{2} ) T + ( 6 m n -
18 n^{4} + 3 n^{2} w^{2} + 2m w \varepsilon ), $$ the discriminant
$ \Delta ( h ) = ( w^{2} - 9 n^{2} )^{2} - 24 m n - 8 m w \varepsilon $ must be a
non-zero square integer. So $ ( w^{2} - 9 n^{2} )^{2} - 24 m n - 8 m w \varepsilon = v ^{2} $
for some $ v \in \Z \setminus \{ 0 \}.  $ \ We denote \ $ t =
w^{2} - 9 n^{2}, $ \ then $$ t^{2} - 24 m n - 8 m w \varepsilon =
v^{2}.  \eqno(2.6) $$ Obviously, \ $ t \equiv v \ ( \text {mod }
\ 2 ). $ \ Let \ $ t = 2 u + v $ \ with \ $ u \in \Z. $
Substituting into equation (2.6), we get
$$ \left \{
   \begin{array}{l}
  u^{2} + uv - 6 m n - 2 m w \varepsilon = 0 , \\
  w^{2} = 9 n^{2} + 2 u + v.
  \end{array}
  \right.   \eqno(2.7) $$
As $ x $ and $ x_{2} $ are the roots of $ h ( T ) , $ we have $$
 x, \ x_{2} = \frac{ ( w^{2} - 3 n^{2} ) \pm \sqrt{ \Delta ( h ) }}{2}
 = \frac{ 6 n^{2} + 2 u + v \pm v }{2} = 3 n^{2} + u \ \text { or } \
 3 n^{2} + u + v . $$ Without loss of generality, we may assume
 that $$ x = 3 n^{2} + u  \quad \text{ and } \quad  x_{2} =  3 n^{2} + u + v. $$
On the other hand, by the duplication formula (2.1) for \ $ 2
P_{9}, $ \ we have  \\
$ x^{4} -2a x^{2} - 8 b x + a^{2} = 4 x_{2} ( x^{3} + ax + b ) =
4 ( w^{2} - 3 n^{2} - x ) ( x^{3} + ax + b ), $  that is
$$ 5 x^{4} - 4 ( w^{2} - 3 n^{2}) x^{3} + 2 a x^{2} - 4 (a ( w^{2} - 3 n^{2}) + b) x
+ a^{2} - 4 b ( w^{2} - 3 n^{2} ) = 0.  \eqno(2.8) $$
Substituting $ x = 3 n^{2} + u, \ a = 6 m n - 27 n^{4}, \ b =
m^{2} - 18 m n^{3} + 54 n^{6} \ $ into the equation (2.8), by a
tedious calculation, we obtain that
$$ \begin{array}{l}
5 u^{4} + ( 72 n^{2} - 4 w^{2} ) u^{3} +
( 324 n^{4} - 36 n^{2} w^{2} + 12 m n ) u^{2} \\
+ ( 216 m n^{3} - 4 m^{2} - 24 m n w^{2} ) u + 36 m^{2} n^{2} - 4
m^{2} w^{2} = 0.
\end{array}  \eqno(2.9) $$ Then putting \ $ w^{2} = 9 n^{2} + 2 u + v
$ \ into (2.9), we get
$$ \begin{array}{l}
5 u^{4} + ( 36 n^{2} - 8 u - 4 v ) u^{3} +
( 12 m n - 72 n^{2} u - 36 n^{2} v ) u^{2} \\
- ( 4 m^{2} + 24 m n v + 48 m n u ) u - 8 m^{2} u - 4 m^{2} v = 0.
\end{array}  \eqno(2.10) $$ Furthermore, by (2.7),
$$ \begin{array}{l}
8 m^{2} u + 4 m^{2} v = 4m^{2} ( 2 u + v ) = 4 m^{2} w^{2} - 36
m^{2} n^{2} \\
= ( 2 m w \varepsilon )^{2} - 36 m^{2} n^{2} = ( u^{2} + u v - 6 m
n )^{2} - 36 m^{2} n^{2} \\
= ( u^{2} + u v )^{2} - 12 m n ( u^{2} + u v ).
\end{array}   \eqno(2.11)  $$ Substituting (2.11) into (2.10), we get
$$ \begin{array}{l}
5 u^{4} + ( 36 n^{2} - 8 u - 4 v ) u^{3} +
( 12 m n - 72 n^{2} u - 36 n^{2} v ) u^{2} \\
- ( 4 m^{2} + 24 m n v + 48 m n u ) u - ( u^{2} + u v )^{2} + 12 m
n ( u^{2} + u v ) = 0.
\end{array}   $$  So \ $ u = 0 $ \ or
$$ \begin{array}{l}
5 u^{3} + ( 36 n^{2} - 8 u - 4 v ) u^{2} +
( 12 m n - 72 n^{2} u - 36 n^{2} v ) u \\
- ( 4 m^{2} + 24 m n v + 48 m n u ) - u ( u +  v )^{2} + 12 m n (
u + v ) = 0.
\end{array}   \eqno(2.12) $$  If \ $ u = 0 , $ \ then by (2.7), \
$ w \varepsilon = -3n $ \ since $ m \neq 0 . $ And then $ v =
w^{2} - 9 n^{2} - 2 u = 0 . $ \ A contradiction! Therefore \ $ u
\neq 0 $ \ and we have the equality (2.12). After expanding and collecting terms,
(2.12) reduces to $$ u \left[ ( 2 u + v )^{2} + 2uv  \right] + 4 ( m + 3 n u
) \left[ m + 3 n ( u + v ) \right] = 0. \eqno(2.13) $$
Furthermore, if necessarily,  replacing \ $ m, \ n $ \ by \ $ m
\varepsilon , \ n \varepsilon $ \  respectively, the symbol \ $
\varepsilon $ \ can be omitted in our equations. Hence all the
conditions in Case (II) are obtained.   \\
Conversely, if the conditions on $ a, b $ in (II) hold, then
from the above procedure, it is easy to verify that the point \ $
 P_{9} = ( 3 n^{2} + u, \ u w - m ) \in E ( \Q )_{\text{tors}} $
 is of order 9. In fact, $ 9 P_{9} = O $ but $ 3 P_{9} = (3 n^{2}, m) \neq O, $
so $ P_{9} $ has exact order $ 9. $ Also
$ 2 P_{9} = ( 3 n^{2} + u + v, \ ( u + v ) w - m ), $
and the case (II) is proved.
\par \vskip 0.2cm

(III) \ If \ $ E(\Q)_{\text{tors}} = \Z / 5 \Z. $ \ Then there is a
point \ $ P_{5} =( x, \ y ) \in E ( \Q )_{\text{tors}}  $ \ such that \ $
5 P_{5} = O $ \ and \ $ P_{5} \neq O. $ \ So \ $ 4P_{5} = - P_{5}
, $ \ and then \ $ x ( 4P_{5} ) = x ( - P_{5} ) = x( P_{5} ) = x.
$ \ Denote \ $ 2 P_{5} = ( x_{2} , \ y_{2} ). $ \ Obviously \ $
x_{2} \neq x. $ \ Applying the duplication formula (2.1) for \ $
4 P_{5}, $  $$ x ( 4 P_{5} ) = \frac{ x_{2}^{4} - 2a x_{2}^{2} -
8 b x_{2} + a^{2}}{ 4 x_{2}^{3} + 4 a x_{2} + 4 b }.  $$ Since \
$ x ( 4P_{5} ) = x, $ \ we get $$ x_{2}^{4} - 2a x_{2}^{2} - 8 b
x_{2} + a^{2} = 4 x x_{2}^{3} + 4 a x x_{2} + 4 b x.  \eqno(2.14)
$$ Also by (2.1) for \ $ 2 P_{5}, $ \ we have
$$ x^{4} - 2a x^{2} - 8 b x + a^{2} = 4 x^{3} x_{2} + 4 a x x_{2} + 4 b x_{2} .
 \eqno(2.15) $$  (2.14) - (2.15), and dividing by \ $ x - x_{2}, $ \
 we obtain $$ x^{3} - 3 x^{2} x_{2} - 3 x x_{2}^{2} + x_{2}^{3} -
 2 a ( x + x_{2} ) - 4 b = 0, \quad \text{that is,} $$
 $$ ( x + x_{2} )^{3} - 6 ( x + x_{2} ) x x _{2} - 2 a ( x + x_{2} ) - 4 b = 0 .
  \eqno(2.16) $$ So \ $ 2 \mid ( x + x_{2} ) . $ \ Let $$ x + x_{2} = 2 u,
  \quad  \quad  x x_{2} = t , \quad  \quad u, \ t \in \Z.  \eqno(2.17) $$
Substituting into (2.16), then follows $$ 2 u^{3} - 3 u t - a u = b. \eqno(2.18) $$
Obviously \ $ u \neq 0. $ \ Otherwise \ $ b = 0, $ \ and then \ $ ( 0, \ 0 ) \in E( \Q ) [2],  $
contradicts to our assumption! So \ $ u \mid b. $ \ Let \ $ b = u r $ \ with $ r
\in \Z \setminus \{0\}. $ \ Then from (2.18) we have \ $ a = 2
u^{2} - ( 3 t + r ). $ \ Since $ x $ and $ x_{2} $ are the two
distinct integer roots of the polynomial \ $ h( T ) = T^{2} - 2 u
T + t , $ \ the discriminant \ $ \Delta ( h ) = ( - 2 u )^{2} - 4
t = 4 ( u^{2} - t ) $ \ must be a non-zero square integer, so \ $
u^{2} - t = v^{2} $ \ for some \ $ v \in \Z \setminus \{0\}. $
\ And then $$ x, \ x_{2} = \frac{ 2 u \pm \sqrt{ \Delta ( h )}}{2}
= \frac{ 2 u \pm 2 v }{2} = u \pm v. $$ Without loss of
generality, we may assume that $$ x = u + v , \quad \text{ and }
\quad x_{2} = u - v.  \eqno(2.19) $$ Now \ $ t = x x_{2} = u^{2} -
v^{2}, $ \ by definition, \ $ y^{2} = x^{3} + a x + b = ( u + v
)^{3} + ( 2 u^{2} - ( 3 t + r ) ) ( u + v ) + u r = 4 v^{3} + 6 u
v^{2} + 2 u^{2} v - r v. $ \ So \ $ v \mid y^{2} . $ \ Let \ $
y^{2} = v s $ \ with \ $ s \in  \Z . $ \ Then \ $ s \neq 0, $ \
and \ $ vs = 4 v^{3} + 6 u v^{2} + 2 u^{2} v - r v . $ \ So \ $ r
+ s = 4 v^{2} + 6 u v + 2 u^{2} = 2 ( u + v ) ( u + 2 v ) . $ \
Denote \ $ y = e \in \Z \setminus \{ 0 \}. $ \ Then \ $ e^{2} = v s. $ \ So
$$ \left \{
   \begin{array}{l}
  a = 2 u^{2} - ( 3 ( u^{2} - v^{2} ) + r ) = 3 v^{2} - u^{2} - r  \\
  \hspace{0.3cm}= 3 v^{2} - u^{2} - ( 2 ( u + v ) ( u + 2 v ) - s
  ) \\
  \hspace{0.3cm}= - v^{2} - 3 u^{2} - 6 u v + s , \\
  b = u r = u ( 4 v^{2} + 6 u v + 2 u^{2} - s ).
  \end{array}
  \right.   \eqno(2.20) $$
Note that \ $ x_{2} = u - v, $ \ by (2.14) we have
$$    \begin{array}{l}
  ( u - v )^{4} - 2 ( - v^{2} - 3 u^{2} - 6 u v + s ) ( u - v )^{2}
   - 8 u ( 4 v^{2} + 6 u v + 2 u^{2} - s ) ( u - v ) \\
  + ( - v^{2} - 3 u^{2} - 6 u v + s )^{2} = 4 ( u + v ) ( u - v )^{3}
  +  \\
  4 ( - v^{2} - 3 u^{2} - 6 u v + s ) ( u + v ) ( u - v ) +
  4 u ( 4 v^{2} + 6 u v + 2 u^{2} - s ) ( u + v ).
  \end{array}  \eqno(2.21) $$
  Via a complicated calculation, (2.21) can be simplified to be
  $$ s^{2} - 12 u  v s + 4 v^{4} = 0.  \eqno(2.22) $$
View (2.22) as a quadratic equation in variable \ $ s. $ \ Since \
$ s \in \Z, $ \ the discriminant $$ \Delta _{s} = ( - 12 u v
)^{2} - 4 \cdot 4 v^{4} = 16 v^{2} ( 9 u^{2} - v^{2}) $$ must be
a square integer. So \ $ 9 u^{2} - v^{2} = k^{2} $ \ for some \ $
k \in \Z . $ \ Then $$ s = \frac{12 u v \pm \sqrt{ \Delta
_{s}}}{2} = 6 u v \pm 2 k v . $$ So we may take \ $ s = 6 u v + 2
k v = 2v ( 3 u + k ). $ \ Then \ $ e^{2} = v s = 2 v^{2} ( 3 u +
k ).  $ \ So \ $ 2 ( 3 u + k ) = ( e / v )^{2} \in 2 \Z , $ and so
\ $ e / v = 2 w $ \ for some \ $ w \in \Z \setminus \{ 0 \}. $ \
Thus \ $ s = v \cdot ( e / v )^{2} = v \cdot 4 w^{2}, $ \ and \ $
e = 2 v w. $ \ Substituting them into (2.22), and note that \ $ v
\neq 0, $ \ we get
$$ v^{2} - 12 u w^{2} + 4 w^{4} = 0. \eqno(2.23) $$ Since \ $ v w \neq 0,
$ \ from (2.23) we have \ $ 4 w^{2} \mid v^{2}, $ \ so \ $ 2w
\mid v . $ \ Let \ $ v = 2 w \theta $ \ with \ $ \theta \in \Z
\setminus \{ 0 \}. $ \ Then by (2.23) we get $$ 3 u = w^{2} +
\theta ^{2}.  \eqno(2.24) $$ From the equality (2.24), we have $ w^{2} + \theta^{2}
\equiv 0 (\text{mod} 3), $ which deduces $  w \equiv \theta
\equiv 0 (\text{mod} 3). $ So \ $ \theta = 3 m $ \ and \ $ w = 3 n , $
where \ $ m,  n \in \Z $ \ and \ $ m n \neq 0. $ Therefore
$$ \left \{
   \begin{array}{l}
  u = \frac{1}{3} ( \theta ^{2} + w^{2} ) = 3 ( m^{2} + n^{2} ),  \\
  v = 2 w \theta = 18 m n , \\
  w = 3 n.
  \end{array}
  \right.   \eqno(2.25) $$
And then
$$ \left \{
   \begin{array}{l}
  a = - v^{2} - 3 u^{2} - 6 u v + s = - v^{2} - 3 u^{2} - 6 u v +
  4 v w^{2}  \\
  \hspace{0.3cm}= - ( 18 m n )^{2} - 18 ( m^{2} + n^{2} ) \cdot 18 m n
  - 3 ( 3 ( m^{2} + n^{2} ) )^{2} + 4 \cdot 18 m n \cdot ( 3 n )^{2} \\
  \hspace{0.3cm}= -27 ( m^{4} + n^{4} + 14 m^{2} n^{2} + 12 m^{3} n -
  12 m n^{3} ) \\
  \hspace{0.3cm}= -27 [( m^{2}-n^{2} + 6mn)^{2} - 20 m^{2} n^{2}
  ], \\
  b = u ( 4 v^{2} + 6 u v + 2 u^{2} - s ) = u ( 4 v^{2} + 6 u v + 2 u^{2} -
  4 v w^{2} ) \\
  \hspace{0.3cm}= 3 ( m^{2} + n^{2} ) \left[ 4 ( 18 m n )^{2} +
  18^{2} m n (  m^{2} + n^{2} ) +
   2 ( 3 ( m^{2} + n^{2} ) )^{2}
  - 72 m n ( 3 n )^{2} \right] \\
  \hspace{0.3cm}= 54 ( m^{2}+ n^{2} ) ( m^{4} + n^{4} + 74 m^{2} n^{2} + 18 m^{3} n -
  18 m n^{3} ) \\
  \hspace{0.3cm}= 54 ( m^{2}+ n^{2} ) \left[ ( m^{2}-n^{2} + 9mn )^{2} - 5 m^{2} n^{2} \right]
  \end{array}
  \right.   \eqno(2.26) $$
  Moreover, $$ P_{5} = ( u + v, \ e ) = ( u + v, \ 2 vw ) =
  ( 3 ( m^{2} + n^{2} ) + 18 m n, \ 108 m n^{2} ), $$
  $$ 2 P_{5} = (u - v, v^{2} / w) =
  (3 ( m^{2} + n^{2}) - 18 m n, \ 108 m^{2} n). $$ So we
  obtain all the conditions in case (III).  \\
Conversely, if the conditions in case (III) hold, then it is easy to verify that the point
\ $ P_{5} = (3 ( m^{2} + n^{2} ) + 18 m n, 108 m n^{2} ) \in E(\Q)_{\text{tors}} $
is of order 5. This proves case (III).
\par  \vskip 0.2 cm

( IV ) \ If \ $ E(\Q)_{\text{tors}} = \Z / 7 \Z. $ \ Then there is a
point \ $ P_{7} =( x, y) \in E(\Q)_{\text{tors}} $ \ such that \ $
7 P_{7} = O $ \ and \ $ P_{7} \neq O. $ \ So \ $ 4P_{7} = - 3
P_{7} , $ \ and then \ $ x ( 4 P_{7} ) = x ( - 3 P_{7} ) = x( 3
P_{7} ). $ \ Denote \ $ 3 P_{7} = ( x_{3}, y_{3} ), 2 P_{7}
= ( x_{2}, y_{2}). $ \ Obviously \ $ 2 P_{7} \neq \pm P_{7} $
\ i. e., \ $ x_{2} \neq x. $  So By the addition law (see [S],
pp.58$ \sim $59),
$$ x_{3} = x ( P_{7} + 2 P_{7} ) = \left(
\frac{y_{2} - y }{ x_{2} - x } \right)^{2}
 - x - x_{2}.  $$
 Let $ \frac{y_{2} - y }{ x_{2} - x } = t,  $
 \ then \ $ t \in \Z $ because $ x, x_{2}, x_{3}, y, y_{2} \in \Z. $
 \ Hence $$ x + x_{2} + x_{3} = t^{2} \quad \text{ and } \quad y_{2}
 - y = t ( x_{2} - x ) .  \eqno(2.27) $$ Moreover, by our
 assumption, it is obvious that $$ x_{3} \neq x, \quad \quad  x_{3} \neq x_{2}.
  \eqno(2.28) $$
  Applying the duplication formula (2.1) for \ $
4 P_{7} , $   $$ x ( 4 P_{7} ) = \frac{ x_{2}^{4} - 2a x_{2}^{2} -
8 b x_{2} + a^{2}}{ 4 x_{2}^{3} + 4 a x_{2} + 4 b }.  $$ Since \
$ x ( 4 P_{7} ) = x (3 P_{7}) = x_{3}, $ \ we get $$ x_{2}^{4} -
2a x_{2}^{2} - 8 b x_{2} + a^{2} = 4 x_{3}( x_{2}^{3} +  a x_{2} +
b ) = 4 ( t^{2} - x - x_{2} ) ( x_{2}^{3} +  a x_{2} + b  ).
\eqno(2.29) $$ Also, by the duplication formula (2.1) for \ $ 2
P_{7}, $ \ we get $$ x^{4} - 2a x^{2} - 8 b x + a^{2} = 4 x^{3}
x_{2} + 4 a x x_{2} + 4 b x_{2}. \eqno(2.30) $$ By definition,
$$ y^{2} = x^{3} + a x + b,  \quad \quad y_{2}^{2} = x_{2}^{3} + a x_{2} + b, $$
so \ $ y_{2}^{2} - y^{2} = x_{2}^{3} - x^{3} + a ( x_{2} - x ). $
\ Since \ $ y_{2} - y = t ( x_{2} - x ) $ \ and \ $ x_{2} \neq x,
$ \ we obtain $$ t ( y + y_{2} ) = x^{2} + x x_{2} + x_{2}^{2} + a. \eqno(2.31) $$ So
$$ \left \{
   \begin{array}{l}
2ty = x^{2} + x x_{2} + x_{2}^{2} + a + t^{2} ( x - x_{2} ) , \\
2ty_{2} = x^{2} + x x_{2} + x_{2}^{2} + a - t^{2} ( x - x_{2} ).
\end{array}
  \right.   \eqno(2.32) $$
From (2.29) and (2.30), it is easy to see that \ $ x \equiv a \ (\text{mod} 2 ) $ and
$ x_{2} \equiv a ( \text{mod} 2), $ \ so \ $ x \equiv x_{2} \ (\text{mod} 2 ). $ \ Hence
we may write \ $ x + x_{2} = 2 \lambda , \ x x_{2} = \gamma $
\ with \ $ \lambda , \ \gamma  \in \Z. $ \ Since $ x $ and $ x_{2}
$ are the two distinct integer roots of the polynomial \ $ h( T )
= T^{2} - 2 \lambda T + \gamma , $ \ the discriminant \ $ \Delta (
h ) = ( - 2 \lambda )^{2} - 4 \gamma = 4 ( \lambda ^{2} - \gamma )
$ \ must be a non-zero square integer, so \ $ \lambda ^{2} -
\gamma = s^{2} $ \ for some \ $ s \in \Z \setminus \{ 0 \}. $ \
Hence
$$ x, \ x_{2} = \frac{ 2 \lambda \pm \sqrt{ \Delta ( h )}}{2} = \frac{ 2
\lambda \pm 2 s }{2} = \lambda \pm s . $$ Without loss of
generality, we may assume that $$ x = \lambda + s , \quad \text{
and } \quad x_{2} = \lambda - s .  \quad  ( s \neq 0 ) \eqno(2.33)
$$ Substituting them into (2.32), we get \ $ 2 t y = 3 \lambda ^{2} +
s^{2} + 2 s t^{2} + a . $ \ Denote \ $ y = e \in \Z \setminus \{ 0
\}. $ \ Then
$$ \left \{
   \begin{array}{l}
 a = 2 e t - 3 \lambda ^{2} - s^{2}  - 2 s t^{2}, \\
 b = y^{2} - x^{3} - a x = e^{2} - ( \lambda + s )^{3} -
 ( 2 e t - 3 \lambda ^{2} - s^{2}  - 2 s t^{2})( \lambda + s ) \\
  \hspace{0.3cm} = e^{2} + 2 \left[ \lambda ^{3} +
  ( s t^{2} -e t - s^{2} ) \lambda  - e s t + s^{2} t^{2} \right],
\end{array}
  \right.   \eqno(2.34) $$ where \ $ \lambda , e, s, t \in \Z
  $ \ and \ $ es \neq 0. $ \ From (2.27), \ $ y_{2} = y + t ( x_{2} - x )
   = e - 2 s t. $ \ So $$ P_{7} = ( \lambda + s, e ) \quad  \text{ and }
  \quad 2P_{7} = ( \lambda - s, e - 2 s t).  \eqno(2.35) $$
  Substituting them into (2.29) and (2.30), respectively, we obtain
$$   \begin{array}{l}
   5 ( \lambda - s )^{4} + 2 ( 2 e t - 3 \lambda ^{2} - s^{2}  - 2 s t^{2} )
    ( \lambda - s )^{2} - \\
    4 \left[ e^{2} + 2 \left[ \lambda ^{3} +
  ( s t^{2} -e t - s^{2} ) \lambda  - e s t + s^{2} t^{2} \right]
  \right] ( \lambda - s ) + \\
  ( 2 e t - 3 \lambda ^{2} - s^{2}  - 2 s t^{2} )^{2} = 4 ( e - 2 s t )^{2}
  ( t^{2} - s - \lambda ),
\end{array}
 \eqno(2.36) $$ and
$$   \begin{array}{l}
   ( \lambda + s )^{4} - 2 ( 2 e t - 3 \lambda ^{2} - s^{2}  - 2 s t^{2}
   ) ( \lambda + s )^{2} - \\
   8 \left[ e^{2} + 2 \left[ \lambda ^{3} +
  ( s t^{2} -e t - s^{2} ) \lambda  - e s t + s^{2} t^{2} \right]
  \right] ( \lambda + s ) + \\
  ( 2 e t - 3 \lambda ^{2} - s^{2}  - 2 s t^{2} )^{2} = 4 e^{2}
  ( \lambda - s  ).
\end{array}     \eqno(2.37) $$
By a tedious calculation, (2.36) and (2.37) can be
simplified to be the following two equalities respectively:
$$   \begin{array}{l}
9 s^{2 }\lambda ^{2} + 6 ( s^{2} t^{2} - s^{3} - e s t ) \lambda +
s^{4} + 6 s^{3} t^{2} \\
- s^{2} ( 3 t^{4} + 6 e t ) + ( 2 e^{2} + 2 e t^{3} ) s = 0 ,
\end{array}   \eqno(2.38) $$
$$   \begin{array}{l}
9 s^{2 }\lambda ^{2} + ( 6 s^{3} - 6 s^{2} t^{2} + 6 e s t - 3
e^{2} ) \lambda + s^{4} - 2 s^{3} t^{2} \\
+  ( t^{4} + 2 e t ) s^{2} - ( e^{2} + 2 e t^{3} ) s + e^{2} t^{2}
= 0.
\end{array}  \eqno(2.39) $$
Via (2.38) - (2.39), we get
$$   \begin{array}{l}
( -12 s^{3} + 12 s^{2} t^{2} - 12 e s t + 3 e^{2}) \lambda + 8 s^{3} t^{2} \\
- ( 4 t^{4} + 8 e t ) s^{2} + ( 3 e^{2} + 4 e t^{3} ) s - e^{2}
t^{2} = 0.
\end{array}  \eqno(2.40) $$
Furthermore, by a complicated calculation, (2.39) can be
decomposed into factors as $$ \left[ ( 3 \lambda + s ) - t^{2}
\right] \cdot \left[ ( 3 \lambda + s ) s^{2} - ( e - s t
 )^{2} \right] = 0 .  \eqno(2.41) $$ So  $$ 3 \lambda + s = t^{2}
 \quad \text{or} \quad  ( 3 \lambda + s ) s^{2} = ( e - s t
 )^{2}.  \eqno(2.42) $$
 If $ 3 \lambda + s = t^{2}, $ \ then by (2.27) and (2.33),
 we have \ $ x_{3} = t^{2} - x - x_{2} = (3 \lambda + s ) - 2 \lambda =
 \lambda + s = x, $ \ contradicts to (2.28). Therefore we must
 have  $ ( 3 \lambda + s ) s^{2} = ( e - s t )^{2}. $
 Then $ 3 \lambda + s = t_{1}^{2} $ for some \ $ t_{1} \in \Z. $ And then
 $ s t_{1} = e - s t. $ So \ $ e = s ( t + t_{1} ).
 $ \ Since \ $ e s \neq 0, $ \ we have \ $ t + t_{1} \neq 0. $ \
 Substituting \ $ e $ \ and \ $ 3 \lambda = t_{1}^{2} - s $ \ into
 (2.40), note that \ $ s \neq 0, $ \ by simplifying we get
 $$ 4 s^{2} + 2 ( t^{2} - t_{1}^{2} ) s + ( t_{1}^{4} - t^{4} + 2 t^{3} t_{1}
  - 2 t t_{1}^{3} ) = 0 , \quad  \text{ i.e., } $$
  $$ 4 s^{2} + 2 ( t + t_{1} ) ( t - t_{1} ) s - ( t + t_{1})
  (t - t_{1} )^{3}  = 0.  \eqno(2.43)  $$
  Obviously, \ $ t - t_{1} \neq 0 $ \ and
  \ $ t \equiv t_{1} \ ( \text{mod} \ 2). $ \ Let \ $ t + t_{1} =
  2 u, \quad t - t_{1} = 2 v , \quad u, \ v \in \Z $ \ and \ $ u v \neq 0. $
  \ Then \ $ t = u + v, \ t_{1} = u - v, $ \ and from (2.43) we get
  $$ s^{2} + 2 u v s - 4 u v^{3} = 0.  $$ Thus \ $ 2 \mid s, $ \
  and so \ $ s = 2 w $ \ for some
  $ w \in \Z \setminus \{ 0 \}. $ \ Then $ w^{2} + u v w - u v^{3} = 0. $
  Moreover, \ $ 3 \lambda + 2 w = ( u - v )^{2} , \quad e = s ( t + t_{1} ) =
  2 w \cdot 2u = 4 u w. $ \ So from (2.34),
  $$ \left \{
   \begin{array}{l}
 a = 2 e t - 3 \lambda ^{2} - s^{2}  - 2 s t^{2} \\
  \hspace{0.3cm} = -3 \lambda ^{2} + 4 w ( u^{2} - v^{2} - w ) ,
  \\
 b = e^{2} + 2 \left[ \lambda ^{3} +
  ( s t^{2} -e t - s^{2} ) \lambda  - e s t + s^{2} t^{2} \right]
  \\
  \hspace{0.3cm} = 8w^{2} ( u^{2} + v^{2} - \lambda ) +
  2 \lambda ( \lambda ^{2} + 2w ( v^{2} - u^{2} )).
\end{array}
  \right.   $$ Also
  $ P_{7} = (\lambda + s, e) = (\lambda + 2 w,  4 u w), \
  2 P_{7} = (\lambda - s, e - 2 s t) = (\lambda - 2 w, -4 v w). $
  Therefore, we obtain all the conditions in case (IV). \\
  Conversely, if the conditions in case (IV) hold, then it is not
  difficult to verify that the point \ $ P_{7} = ( \lambda + 2 w , \ 4 u w )
  \in E(\Q)_{\text{tors}} $ is of order 7. So case (IV) is proved.
  This completes the proof of Theorem 1.
 \par  \vskip 0.2 cm

{\bf Remark.} \ This paper is a revised version of the early one ([Q2], 2008).
\par  \vskip 0.2 cm

{ \bf Acknowledgments } \ We would like to thank the referee
for careful reading and helpful suggestions.

\par  \vskip 0.6 cm

\hspace{-0.8cm} {\bf  References }
\begin{description}

\item[[F1]] Y. Fujita, Torsion subgroups of elliptic curves with
non-cyclic torsion over $ \Q $ in elementary abelian 2- extensions
of $ \Q , $ { \it Acta Arithmetica}, 2004, 115: 29-45.

\item[[F2]] Y. Fujita, Torsion subgroups of elliptic curves in
elementary abelian 2- extensions of $ \Q , $ {\it Journal of
Number Theory}, 2005, 114: 124-134.

\item[[F3]] Y. Fujita, The 2-primary torsion on elliptic curves in
the $ \Z_{p}- $ extensions of $ \Q , $ Manuscripta mathematica,
2005, 118:339-360.

\item[[Kn]] A. Knapp, Elliptic Curves, Princeton: Princeton Univ.
Press, 1992.

\item[[Kw]] S. Kwon, Torsion subgroups of elliptic curves over
quadratic extensions, {\it Journal of Number Theory}, 1997, 62:
144-162.

\item[[M1]] B. Mazur, Modular curves and the Eisenstein ideal,
 {\it IHES  Publ. Math}., 1977, {\bf 47}: 33-186.

\item[[M2]] B. Mazur, Rational points on modular curves, Modular
Functions of One Variable V, {\it Lecture Notes in Math}., New
York: Springer-Verlag, 1977, {\bf 601}: 107-148.

\item[[O]] K. Ono, Euler's concordant forms, {\it Acta
Arithmetica}, 1996,  {\bf LXX VIII}(2): 101-123.

\item[[Q1]] D.R. Qiu, Mordell-Weil groups and related problems of
elliptic curves, { \it PhD Thesis, Tsinghua University}, 2000.

\item[[Q2]] D.R. Qiu, An explicit classification for the cyclic rational torsion subgroups
of odd order of elliptic curves over $ \Q. $ arXiv: 0803.0058 v1, 2008.

\item[[QZ1]] D.R. Qiu, X.K. Zhang, Explicit classification for
torsion cyclic subgroups of rational points with even orders of
elliptic curves, { \it Chinese Science Bulletin}, 1999,
44(21):1951-1952.

\item[[QZ2]] D.R. Qiu, X.K. Zhang, Explicit classification for
torsion subgroups of rational points of elliptic curves, { \it
Acta Mathematica Sinica, English Series}, 2002, 18(3):539-548.

\item[[QZ3]] D.R. Qiu, X.K. Zhang, Elliptic curves and their
torsion subgroups over number fields of type (2, 2, ..., 2), { \it
Science in China ( Series A)}, 2001, 44(2): 159-167.

\item[[S]] J. Silverman, The Arithmetic of Elliptic Curves, New
York: Springer-Verlag, 1986.

\end{description}

\end{document}